\documentstyle{amsppt}
\magnification=1200
\NoBlackBoxes
\hcorrection{0.25in}   %?
\TagsOnRight
\nologo

\def\ref#1{[#1]}
\define\wt{\widetilde}

\define\ra{\rightarrow}
\define\bs{\backslash}
\define\ov{\overline}
\define\col{\!:\!}

\define\br{{\Bbb R}}
\define\bz{{\Bbb Z}}

   %%%%

\define\Hom{\text{Hom}}
\define\fix{\text{fix}}

\define\isom{\operatorname{Isom}}

\define\Ga{\Gamma}

\define\La{\Lambda}
\define\p{\partial}
\define\a{\alpha}

\define\Sa{\Sigma}

\define\la{\lambda}
\define\om{\Omega}

\define\sch{{\Cal H}}

\define\sci{{\Cal I}}

\define\bb{{\Bbb B}}
\define\bc{{\Bbb C}}

\define\bh{{\Bbb H}}

\define\ch{\bh_{\bc}^}
\define\cb{\Bbb B_{\bc}^n}

\define\im{\operatorname{Im}}

\parindent=1pc                  %% AMSPPT default setting

% added by arXiv admin 2009-09-07
\voffset=-0.8in

\leftheadtext{Boris Apanasov }
\rightheadtext{Quasiconformality and Geometrical Finiteness}
\topmatter
\title Quasiconformality and Geometrical Finiteness in
Carnot--Carath\'eodory and Negatively Curved Spaces\endtitle

\author Boris Apanasov\footnote"{}"{Supported in part by the NSF. Research at MSRI is supported in part by NSF grant DMS-9022140.}
\endauthor
\address Department of Mathematics, University of Oklahoma,
Norman, OK  73019
\endaddress
\email apanasov\@uoknor.edu \endemail
\subjclass 57,55,53,51
\endsubjclass
\abstract\nofrills {ABSTRACT.}
The paper sketches a recent progress and formulates several open
problems in studying equivariant quasiconformal and quasisymmetric 
homeomorphisms in negatively curved spaces as well as
 geometry and topology of noncompact geometrically finite 
negatively curved manifolds and their boundaries at infinity having 
Carnot--Carath\'eodory structures.
Especially, the most interesting are complex hyperbolic manifolds with 
Cauchy--Riemannian structure at infinity, which  occupy a distinguished
niche and whose properties make them surprisingly different from real
hyperbolic ones.
 \endabstract
\endtopmatter
\bigskip

\document
\head 1. Introduction \endhead

The paper sketches a recent progress and formulates several open
problems in studying quasiconformal and quasisymmetric homeomorphisms as well as
geometry and topology of noncompact geometrically finite 
negatively curved manifolds and their boundaries at infinity having Carnot--Carath\'eodory structures 
in the sense of M.~Gromov \cite{Gr}. Here, complex hyperbolic and Cauchy--Riemannian 
manifolds occupy a distinguished
niche. First, due to their complex analytic nature, a broad spectrum of 
techniques can contribute to the study, and already obtained results show 
surprising differences between geometry and topology of noncompact complex and real 
hyperbolic manifolds, see \cite{BS, BuM, EMM, GoM, 
KR1-2}. Second, these K\"ahler manifolds are
the most known manifolds with variable pinched negative 
curvature \cite{BGS, B, Mg1, MaG}.
Finally, for complex analytic surfaces, one can apply 
Seiberg-Witten invariants, decomposition of 4-manifolds along
homology 3-spheres, Floer homology and new (homology) cobordism
invariants \cite{LB, A11, Sa1-Sa4}. 

There are three
main themes united by using of general Thurston's idea of geometric
uniformization of low dimensional manifolds provided some canonical
decomposition of them into geometric pieces. First, basic geometric block-manifolds
of negative (variable) curvature which possess geometrical finiteness and may
however have infinite volume. Second, the interaction between the negatively 
curved geometry of such blocks and the induced Carnot--Carath\'eodory structure of their boundary
manifolds at infinity, especially the interaction between K\"ahler geometry of
geometrically finite complex hyperbolic manifolds and Cauchy--Riemannian structure of 
their boundaries modeled on the Heisenberg group (which is a particular case of 
Carnot groups). And finally, deformations of negatively curved 
manifolds inducing automorphisms of boundaries at infinity which 
preserves their natural contact structures, together with their metrical (quasiconformal
or quasisymmetric) properties. 

Our study of geometrical finiteness and homeomorphisms in spaces with negative 
(variable) curvature 
exploits a new structural theorem about isometric actions of
discrete groups on nilpotent Lie groups (in particular on the 
Heisenberg group $\sch_n$). % which generalizes known Bieberbach and Auslander 
%results \cite {Au}. 
The problems there have a unique appeal both for the amount
of similarity with the model situation of interactions between real hyperbolic 
geometrically finite manifolds, their boundaries with natural conformal structures 
and their quasiconformal deformations (see
\cite{A1, A11}), and for the interesting ways in which the similarity 
breaks down. 

Besides geometrical finiteness in negatively curved spaces, one of inspiring ideas 
of our study is a well known theorem of D.Sullivan, which gave rise to many important results
in geometry and topology of manifolds and theory of quasiconformal mappings. It states,
 see \cite{Su1, TV}, that 
homeomorphisms of quasiconformal $n$-manifolds, $n\neq4$,  can be approximated by 
quasiconformal ones, where one of important classes is the class of 
conformal manifolds. Here by the quasiconformality one means the boundedness 
of distortion with respect to the Euclidean metric. This result rises
questions of approximation in the quasiconformal category in another 
negatively curved geometries and the corresponding sub-Riemannian Carnot--Carath\'eodory structures 
which appear at infinity of those geometries as well. 
 We mention here a recent important development by M.Gromov on 
Carnot--Carath\'eodory spaces \cite{Gr}, where he shown that continues mappings can be approximated 
by mappings that are Lipschitz with respect to the Carnot--Carath\'eodory
metric. Another important achievement is a recent result by G.Margulis
and G.Mos\-tow \cite{MM} that quasiconformal mappings of Carnot--Carath\'eodory spaces are 
a.e. differentiable and preserve their contact structures.

Another set of problems is related to the stability theorem proven by 
D.Sullivan \cite{Su3} for planar Kleinian groups (see also
\cite{A11, MaG}). It raises questions on deformations of discrete
isometry groups in spaces of negative (variable) curvature, varieties of their
representations (Teichm\"uller spaces) and their boundaries, and 
quasiconformal/quasisymmetric homeomorphisms induced by
isomorphisms of geometrically finite groups (in conformal
category, it is known as Tukia isomorphism theorems \cite {Tu}).

\head 2. Complex hyperbolic and Heisenberg manifolds  \endhead

The natural class of manifolds where the above problems can be
considered is the class of geometrically finite locally homogeneous
manifolds. That is why we start with studying such manifolds and 
their boundaries at infinity, especially complex hyperbolic manifolds
and Cauchy--Riemannian and Heisenberg manifolds at their infinity.
 
We recall some facts
concerning the link between nilpotent geometry of the Heisenberg group
and the K\"ahler geometry of the complex hyperbolic space (see
\cite{Go, GP1, KR1-2}).  Let $\cb\subset\bc^n$ be the unit complex ball equipped 
with the Bergman metric $d$ which turns the ball into the complex hyperbolic $n$-space
$\ch n$ whose sectional curvature is between -1 and -1/4. The automorphism group
of $\ch n\cong\cb$ is the projective unitary group
$PU(n,1)$ whose elements $g\in PU(n,1)$ are biholomorphisms of $\bb^n_{\bc}$
of the following three types. 
If $g$ fixes a point in $\ch n$, it is called {\it elliptic}.
If $g$ has exactly one fixed point, and it lies in
$\partial \ch n$, $g$ is called {\it parabolic}. If $g$
has exactly two fixed points, and they lie in $\partial \ch n$,
$g$ is called {\it loxodromic}. These three types exhaust all the possibilities.
Since $PU(n,1)$ can be embedded in a linear group (for details, see
\cite{AX1}), any finitely generated group $G\subset PU(n,1)$
is residually finite and has a finite index torsion free subgroup.

For any point $\infty\in\p\cb$, one can identify $\ov{\cb}\bs\{\infty\}$
with the closure of the Siegel domain $\frak S_n$, which is conveniently 
represented in horospherical coordinates as 
$\bc ^{n-1}\times \br\times [0,\infty)$ (\cite {GP1}). Then ``the boundary plane"
$\bc ^{n-1}\times \br\times \{0\}=\p\ch n\bs\{\infty\}$ and the horospheres
$H_u=\bc ^{n-1}\times \br\times \{u\}$, $0<u<\infty$, centered at $\infty$ are
identified with the Heisenberg group $\sch_n=\bc^{n-1}\times \br$. 
It is a 2-step
nilpotent group with center $\{0\}\times \br \subset \bc ^{n-1}\times \br$, and 
the inverse of $(\xi,v)$ is  $(\xi,v)^{-1}=(-\xi,-v)$. 
The Heisenberg group $\sch_n$ isometrically acts
on itself and on $\ch n$ by left translations:

$$T_{(\xi_0,v_0)} \col (\xi,v,u)\longmapsto (\xi_0+\xi\,,v_0+v+2\im
\langle\langle \xi_0,\xi\rangle\rangle \,,u)\,.$$

The unitary group $U(n-1)$  acts on  $\sch_n$
and $\ch n $ by rotations:
$A(\xi,v,u)=(A\xi\,,v\,,u)$ for $A\in U(n-1)$. The semidirect product
$\sch(n)=\sch_n \rtimes U(n-1)$ is naturally embedded in
$U(n,1)$ as follows:

$$A \longmapsto \pmatrix
A&0&0 \\ 0&1&0 \\ 0&0&1 \endpmatrix \in U(n,1)\quad \text{for} \quad A\in U(n-1)\,,$$
$$(\xi,v)\longmapsto \pmatrix
I_{n-1}&\xi & \xi\\
-\bar{\xi}^t&1-{1\over 2}(| \xi | ^2-iv)&
-{1\over 2}(|\xi|^2-iv)\\
\bar\xi^t&{1\over 2}(|\xi|^2-iv)&
1+{1\over 2}(|\xi|^2-iv)\endpmatrix \in U(n,1)$$
where $ (\xi,v)\in \sch_n= \bc ^{n-1}\times \br$ and $\bar \xi^t$ is
the conjugate transpose of $\xi$.

The action of $\sch(n)$ on $\overline{\ch n} \bs \{ \infty \}$
also preserves the Cygan metric $\rho_c$ there, which 
plays the same role as the Euclidean metric does on the upper half-space 
model of the real hyperbolic space $\bh^n$ and is induced by the
following norm:

$$||(\xi,v,u)||_c=|\,||\xi||^2 + u - iv|^{1/2}\,,\quad 
(\xi,v,u)\in\bc ^{n-1}\times \br\times [0,\infty)\,.$$

The relevant geometry on each horosphere $H_u\subset \ch n$, 
$H_u\cong\sch_n=\bc ^{n-1}\times \br$, is 
the spherical $CR$-geometry
induced by the complex hyperbolic structure.
 The geodesic perspective from ${\infty}$
defines $CR$-maps between horospheres, which extend to $CR$-maps
between the one-point compactifications 
$H_u\cup {\infty}\approx S^{2n-1}$.
In the limit, the induced  metrics on horospheres fail to converge 
but the $CR$-structure
remains fixed. In this way, the complex hyperbolic geometry induces
$CR$-geometry on the sphere at infinity $\p \ch n \approx S^{2n-1}$, naturally
identified with the one-point compactification of the Heisenberg group
$\sch_n$.

\head 3. Geometrical finiteness in negative curvature\endhead

Our main assumption on a negatively curved $n$-manifold $M$ is
the geometrical
finiteness condition on its fundamental group $G\subset \isom X$ acting by 
isometries on a simply connected space $X=\wt M$, which in
particular implies (see below) that the discrete group $G$ is finitely generated. 

A subgroup $G\subset \isom X$ is called {\it discrete} if it is a discrete  subset
of $\isom X$. The {\it limit set} $\La(G)\subset \p X$ of a discrete group
$G$ is the set of
accumulation points of (any) orbit $G(y),\, y\in X$. The complement of
$\La(G)$ in $\p X$ is called the {\it discontinuity set} $\om(G)$.
A discrete group $G$ is called {\it elementary}
if its limit set $\La(G)$ consists of at most two points.  
An infinite discrete group $G$ is called {\it parabolic}
if it has exactly one fixed point $\fix (G)$; then $\La(G)=\fix (G)$, and
$G$ consists of either parabolic or elliptic elements. As it was
observed by many authors (c.f. \cite {MaG}), parabolicity in the
variable curvature case is not as easy a condition to deal with as it is
in the constant curvature space. However the results below simplify the
situation.

Due to the absence of totally geodesic 
hypersurfaces in a space $X$ of variable negative curvature,  we cannot use the original
definition of geometrical finiteness which came from
an assumption that the corresponding real hyperbolic manifold $M=\bh^n/G$ may be 
decomposed  into a cell
by cutting along a finite number of its totally geodesic hypersurfaces, that is
the group $G$ should possess a finite-sided fundamental polyhedron, see  
\cite {Ah}.  However, one can use another definition of 
 {\it geometrically finite} groups $G\subset \isom X$
as those ones whose limit sets $\La(G)\subset \p X$ consist of only conical limit points 
and parabolic (cusp) points $p$ with compact quotients $(\La(G)\bs\{p\})/G_p$
with respect to parabolic stabilizers $G_p\subset G$ of $p$, see
\cite {BM, Bow}. There are other definitions of geometrical finiteness 
in terms of ends and the minimal convex retract of the noncompact manifold $M$, 
which work well not only in the real hyperbolic spaces $\bh^n$ 
(see \cite {Md, Th, A3, A1}) but also in spaces with variable pinched negative 
curvature \cite {Bo}. 

In the case of variable curvature, it is problematic to use
geometric methods based on consideration of finite sided
fundamental polyhedra. In particular, Dirichlet
polyhedra $D_y(G)$ in the complex hyperbolic space $X=\ch n$ are 
fundamental polyhedra for a discrete subgroup $G\subset PU(n,1)$.   
They are bounded by bisectors in a complicated way
(see \cite { Mo2, GP1}), and the bisectors are not totally geodesic hypersurfaces. 
For discrete parabolic
groups $G\subset \isom X$, one may expect that the  Dirichlet
polyhedron $D_y(G)$ centered at a point $y$ lying in a $G$-invariant
subspace has finitely many sides. It is true for real hyperbolic spaces
\cite{A1} as well as for cyclic and dihedral parabolic groups in
complex hyperbolic spaces $X=\ch n$. Namely, due to \cite{Ph}, 
Dirichlet polyhedra $D_y(G)$ are always two sided for any cyclic group
$G\subset PU(n,1)$ generated by a Heisenberg translation.
Due to \cite{GP1}, such finiteness also holds
for a cyclic ellipto-parabolic group or a dihedral parabolic group 
$G\subset PU(n,1)$ generated by inversions
in asymptotic complex hyperplanes in $\ch n$ if the central point $y$ lies
in a $G$-invariant vertical line or $\br$-plane (for any other center $y$, 
$D_y(G)$ has infinitely many sides). Our technique easily implies that this
finiteness still holds for generic parabolic cyclic groups \cite{AX1}:

\proclaim{Theorem 3.1}
 For any discrete group $G\subset PU(n,1)$
generated by a parabolic ele\-ment, there exists a point  $y_0\in \ch n$
 such that the Dirichlet polyhedron $D_{y_0}(G)$ centered at $y_0$
has two sides.
 \endproclaim

However, the behavior of Dirichlet polyhedra for parabolic groups 
$G\subset PU(n,1)$ of rank more than one can be very bad. It is given by 
our construction \cite{AX1}:

\proclaim{Theorem 3.2}  Let $G\subset PU(2,1)$ be a discrete parabolic
group conjugate to the following subgroup $\Ga$ of the Heisenberg group
$\sch_2=\bc \times \br$:
$$\Ga=\{(m,n)\in \bc\times \br \col  m,n\in \bz\}\,.$$
Then any Dirichlet polyhedron $D_y(G)$ centered at an arbitrary 
point $y\in \ch 2$ has  infinitely many sides.
\endproclaim

 Despite this, we are providing \cite{AX1}
a construction of fundamental polyhedra $P(G)\subset\ch n$
for arbitrary discrete parabolic groups $G\subset PU(n,1)$, which are
bounded by finitely many
hypersurfaces (different from Dirichlet bisectors). 
This construction may be seen as a base for extension of
Apanasov's construction \cite {A1} of finite sided
pseudo-Dirichlet polyhedra
in $\bh^n$ to the case of the complex hyperbolic space $\ch n$, which
may solve the following problem:

\proclaim{Problem 3.3} Given geometrically finite group $G\subset \isom X$ in a
negatively curved space $X$, is there any finitely sided fundamental
polyhedron $P(G)\subset X$?
\endproclaim

As another tool for attacking this problem, one can use 
the following structural theorem 
for discrete groups acting on nilpotent Lie groups (lying at infinity of $X$), 
in particular on the Heisenberg group $\sch_n, \sch_n \cup \{\infty\}=\p \ch n$ 
(see \cite {AX1- AX3}):

\proclaim{Theorem 3.4} Let $N$ be a connected, simply connected nilpotent Lie group, 
$C$ a compact group of automorphisms of $N$, and $\Ga$ a discrete subgroup of the 
semidirect product $N\rtimes C$. Then there exist a connected Lie subgroup
$V$ of $N$ and a finite index normal subgroup $\Gamma^*$ 
of $\Gamma$ with the following properties:
\roster
\item There exists $b\in N$ such that $b\Gamma b^{-1}$ preserves
$V$.
\item $V/b\Gamma b^{-1}$ is compact.
\item $b\Gamma^* b^{-1}$ acts on $V$ by left translations and this action is free.
\endroster
\endproclaim

Here we remark that 
\roster

\item Compactness of $C$ is an essential condition
because of Margulis \cite {Mg2} construction of nonabelian free 
discrete subgroups $\Gamma$ of $R^3\rtimes GL(3,R)$.

\item This theorem generalizes a result by L.Auslander \cite {Au} who claimed those 
properties not for whole group $\Ga$ but only for some finite index subgroup. 
We also remark that an extension in \cite{AX1} of Wolf's argument \cite {Wo} 
to the complex hyperbolic case (based on  Margulis Lemma \cite {Mg1, BGS} and geometry 
of $\sch_n$) does not work in general nilpotent groups. A different algebraic
proof of Theorem 3.4 can be found in \cite{AX2, AX3}.
\endroster

This result allows one to study 

\proclaim {Problem 3.5} Investigate parabolic cusp ends of negatively
curved manifolds.
\endproclaim

To solve this problem in the complex hyperbolic case, we use the above Theorem 3.4 to 
define standard cusp neighborhoods of parabolic cusps. This allows us to prove the following 
results about cusp ends of  complex hyperbolic manifolds (or, equivalently, about 
the structure of Heisenberg manifolds).

Namely, suppose a point
$p\in \p \ch n$ is fixed by some
parabolic element of the group $G$, and 
$G_p$ is the stabilizer of $p$ in $G$. 
Conjugating the group $G$ by an element $h\in PU(n,1)\,, h(p)=\infty$,
we have that $G_{\infty}\subset
\sch(n)$. In particular, if $p$ is the origin in  $\sch_n$, we can take
such $h$ as the Heisenberg inversion $\sci$ in the hyperchain in
$\sch_n$, which preserves the unit Heisenberg sphere $S_c(0,1)$ and acts in 
$\sch_n$ as follows:

$$\sci(\xi,v)=\left (\frac {\xi}{|\xi|^2-iv}\,,\, \frac {-v}{v^2+|\xi|^4}
\right )\quad \text {where}\,\,(\xi,v)\in \sch_n=\bc^{n-1}\times \br \,.
\tag3.6$$

 Then, due to Theorem 3.4,
there exists a connected Lie subgroup $\sch_\infty
\subseteq \sch_n$ preserved by $G_{\infty}$.

\definition{Definition} A set  $U_{p,r}\subset
\overline{\ch n}\bs \{p\}$ is called a {\it standard cusp neighborhood
of radius} $r>0$ at a parabolic fixed point $p\in \p \ch n$ of a discrete
group $G\subset PU(n,1)$ if, for the Heisenberg inversion 
$\sci_p\in PU(n,1)$ with
respect to the unit sphere $S_c(p,1)=\{(\xi,v)\col \rho_c(p,\,(\xi,v))=1\}$,
$\sci_p(p)=\infty\,,$ the following conditions hold:
\roster
\item $U_{p,r}=\sci_p^{-1}\left(\{x\in \ch n\cup \sch_n\col
\rho_c(x,\sch_\infty)\geq 1/r \}\right)\,;$
\item  $U_{p,r}$ is precisely invariant
with respect to $G_p\subset G$,
that is:
$$\gamma(U_{p,r})=U_{p,r} \quad  \text{for} \,\,\;\gamma\in G_p\quad \text {and}\quad
g(U_{p,r})\cap U_{p,r}=\emptyset \quad \text {for} \,\,\,g\in G\backslash
G_p\,.$$
\endroster
A parabolic point $p\in \p \ch n$ of $G\subset PU(n,1)$ is
a {\it cusp point} if it has a cusp neighborhood $U_{p,r}$.
\enddefinition

We remark that some parabolic points of a discrete group
$G\subset PU(n,1)$ may not be cusp points, see \cite{AX1, \S 5.4}. Applying
Theorem 3.4 and \cite{Bo}, we have:

\proclaim {Lemma 3.7} Let $p\in \p \ch n$ be a parabolic fixed point of
a discrete group $G\subset PU(n,1)$. Then $p$ is a cusp point if and only if
$(\La(G)\bs \{p\})/G_p$ is compact.
\endproclaim

This fact and \cite {Bo} allows us to use another equivalent definitions of geometrical finiteness.
In particular, a group $G\subset PU(n,1)$ is
{\it geometrically finite}  if and only if its quotient space $M(G)=[\ch n\cup\om (G)]/G$ 
has finitely many
ends, and each of them is a cusp end, that is an end whose neighborhood can be
taken as $U_{p,r}/G_p\approx (S_{p,r_0}/G_p) \times (0,1]$, where 

$$S_{p,r}=\sci_p^{-1}\left(\{x\in H_ \bc ^n\cup \sch_n
 \col \rho_c(x,\sch_\infty)=1/r\}\right)\,.$$

Using the above description of discrete group actions on a nilpotent group 
(Theorem 3.4), one can study Carnot--Carath\'eodory manifolds whose fundamental groups act at 
infinity of a negatively  curved space $X$ as discrete parabolic groups. In particular, 
we establish fiber bundle structures on Heisenberg manifolds which
have the form
$\sch_n/G$ where $G$ is a discrete group freely acting 
on $\sch_n$ by isometries, i.e. a torsion free discrete subgroup of
$\sch(n)=\sch_n\rtimes U(n-1)$, see \cite{AX1}: 

\proclaim{Theorem 3.8} Let $\Ga\subset \sch_n\rtimes U(n-1)$ be a torsion-free 
discrete group  acting on the Heisenberg group $\sch_n= \bc ^{n-1}\times \br $ 
with non-compact quotient.
Then the quotient $\sch_n/ \Ga$ has zero Euler characteristic and
is a vector bundle over a  compact manifold. Furthermore,
this compact manifold is finitely covered
by a nil-manifold which is either a torus or the total space of a 
circle bundle over a torus.
\endproclaim

It gives, in addition to finiteness of generators of geometrically finite groups
established in \cite{Bo}, the following important finiteness:

\proclaim{Corollary 3.9} The fundamental groups of Heisenberg manifolds
and geometrically finite complex hyperbolic manifolds are
finitely presented.
\endproclaim

Due to Theorem 3.4, any Heisenberg manifold $N=\sch_n/\Ga$ is
the vector bundle $\sch_n/\Ga\rightarrow \sch_\Ga/\Ga$ where 
$\sch_\Ga\subset \sch_n$ is a minimal $\Ga$-invariant subspace.
 As simple examples in
\cite{AX1} show, such vector bundles are non-trivial in general.
However, up to finite coverings, they are trivial \cite{AX1}:

\proclaim{Theorem 3.10} Let $\Ga\subset \sch_n\rtimes U(n-1)$ be a
discrete group and $\sch_\Ga\subset \sch_n$  a 
connected $\Ga$-invariant Lie subgroup on which $\Ga$ acts co-compactly.
Then there exists a finite index subgroup $\Ga_0\subset \Ga$
such that the vector bundle $\sch_n/\Ga_0\rightarrow \sch_\Ga/\Ga_0$ is
trivial. In particular, any Heisenberg orbifold $\sch_n/\Ga$ is
finitely covered by the product of
a compact nil-manifold $\sch_\Ga/\Ga_0$ and an Euclidean space.
\endproclaim

Such finite covering property holds not only for Heisenberg 
manifolds alone but for geometrically finite complex hyperbolic manifolds, too:

\proclaim{Theorem 3.11} Let $G\subset PU(n,1)$ be a geometrically
finite discrete group.
Then $G$ has a subgroup $G_0$ of finite
index such that every parabolic subgroup of $G_0$ is isomorphic
to a discrete subgroup of  the Heisenberg group $\sch_n=\bc^{n-1}\times \br$ .
In particular, each parabolic subgroup of $G_0$ is free Abelian or
2-step nilpotent.
\endproclaim

It seems this property holds for discrete isometry groups in 
general negatively curved spaces because
our proof of Theorem 3.11 in \cite {AX1} is 
based on the residual finiteness of geometrically 
finite subgroups $G\subset PU(n,1)$
and the following two lemmas (the second of which generalizes a result 
for finite volume real hyperbolic manifolds (\cite { AF}).

\proclaim {Lemma 3.12}  Let $G \subset \sch_n\rtimes U(n-1)$ be a discrete
 group and $\sch_G\subset \sch_n$ a minimal $G$-invariant connected
 Lie subgroup (given by Theorem 4.1).
Then $G$ acts on  $\sch_G$ by translations if $G$ is either
 Abelian or 2-step nilpotent.
\endproclaim

\proclaim{Lemma 3.13} Let $G\subset \sch_n\rtimes U(n-1)$ be a torsion free
discrete group, $F$ a finite group and $\phi\col  G\longrightarrow F$ an
epimorphism.
 Then the rotational part of $ \ker (\phi)$
has strictly smaller order than that of $G$ if one of the following
happens:
 \roster
 \item $G$ contains a  finite index Abelian subgroup and $F$ is not Abelian;
\item $G$ contains a finite index 2-step nilpotent subgroup  and $F$ is not
a 2-step nilpotent group.
\endroster
\endproclaim

 We conclude this section by pointing out that the problem of geometrical
finiteness is very different in complex dimension two.
Namely, it is a well
known fact that any finitely generated
discrete subgroup of $PU(1,1)$
is geometrically finite. This and Goldman's local rigidity theorem
for uniform lattices $G\subset U(1,1)\subset PU(2,1)$ (see \cite {GM})
suggest the following intrigue question:

\proclaim{Problem 3.14} Are all finitely generated discrete groups
 $G\subset PU(2,1)$ with non-empty discontinuity set 
$\om(G)\subset \p\ch 2$ geometrically finite? 
\endproclaim

To solve this problem, one can try to decompose a given 
finitely generated discrete group $G\subset PU(2,1)$ with non-empty 
discontinuity set $\om(G)\subset \p\ch 2$ into free amalgamated products 
of elementary groups. More arguments for an affirmative solution of
this problem are given by the following two facts. First, due to Chen-Greenberg \cite{CG},
all pure loxodromic subgroups $G\subset PU(2,1)$
are discrete. The second fact is due to the trace classification of projective
transformations \cite{Go, VI.2}. Namely, in contrast to Kleinian groups on the plane, 
the subset of groups in a deformation space of a pure loxodromic group 
$G\subset PU(2,1)$ which have accidental parabolic
elements may have real codimension 1.

\head 4. The boundary at infinity of negatively curved manifolds
      \endhead

Another set of problems is related to geometry and topology of 
Carnot--Carath\'eo\-do\-ry 
manifolds which are the boundaries at infinity of negatively curved 
non-compact manifolds. 

In a sharp contrast to the real hyperbolic case, for a compact complex manifold 
$M(G)=(\ch n\cup \om(G))/G$, an application of Kohn-Rossi %\cite { KR} 
analytic extension theorem shows that the boundary of this manifold $M(G)$ 
is connected, and the limit set $\La(G)$ is in some sense small
(see \cite {EMM} and, for quaternionic and Cayley hyperbolic manifolds,
\cite {C1}).
Moreover, according to a recent result of D.Burns (based on an uniformization
theorem \cite { BuM} for isolated ends of complex analytic spaces), the same
claim about connectedness of the boundary $\p M(G)$ still holds if only a
boundary component is compact.

However, if $\p M(G)$ has no compact components, and there is no
finiteness condition on the holonomy group of 
$M(G)$, our algorithmical construction in \cite{AX1} shows that the situation is 
completely different:

\proclaim{Theorem 4.1}  For any integers $k, k_0$, $k\geq k_0\geq 0$,
and $n\geq 2$,
there exists a complex hyperbolic $n$-manifold $M=\ch n/G$,
$G\subset PU(n,1)$, whose boundary at infinity splits up into
$k$ connected $(n-1)$-manifolds,
$\p_\infty M=N_1\cup\cdots\cup N_k$.
Moreover, for each boundary component $N_j$, $j\leq k_0$, the inclusion
$i_j\col N_j\subset M(G)$
induces a homotopy equivalence of $N_j$ to $M(G)$.
\endproclaim
 
The construction in the proof of this theorem firstly provides discrete
groups
$G\subset PU(n,1)$ with two connected components of the discontinuity set 
$\om(G)\subset \p\ch n$. Here we essentially use properties of
Heisenberg inversions (3.6) in hyperchains in $\ch n$. Then we apply
an idea due to A.Tetenov
\cite {Te, KAG} in order to construct groups $G$ with any given
number $k$ of topological 
$G$-invariant balls $\om_i\subset \om(G)$ with common boundary $\p \om_i=\La(G)$. 
The groups we construct in the proof are all however infinitely
generated. The finitely generated case has a close relation to Problem
3.14:

\proclaim {Problem 4.2} How many boundary components are there in a
complex manifold $M(G)$ with finitely generated fundamental group
$G\subset PU(n,1)$?
\endproclaim

Toward this problem,
we can show that the situation described in Theorem 4.1 is
impossible if the complex hyperbolic manifold $M$ is geometrically finite. 
Namely, if the manifold $M(G)$ has non-compact boundary $\p M=\om(G)/G$ with
a component $N_0\subset \p M$ homotopy equivalent
to $M(G)$, then there exists a
compact homology cobordism $M_c\subset M(G)$ homotopy equivalent to $M(G)$, and
$M(G)$ can be easily reconstructed from $M_c$ by gluing up a finite number of
standard  open ``Heisenberg collars", see \cite {AX1}:

\proclaim {Theorem 4.3} Let $G\subset PU(n,1)$ be a geometrically
finite non-elementary torsion free discrete group whose Kleinian
manifold $M(G)$ has non-compact boundary $\p M=\om(G)/G$ with
a component $N_0\subset \p M$ homotopy equivalent
to $M(G)$. Then there exists a
compact homology cobordism $M_c\subset M(G)$ such that $M(G)$ 
can be reconstructed from $M_c$ by gluing up a finite number of
open collars $M_i\times [0,\infty)$ where each $M_i$ is finitely covered by
the product $E_k\times B^{2n-1-k}$ of a closed $(2n-1-k)$-ball and a
closed $k$-manifold $E_k$ which is either flat or a 
nil-manifold (with 2-step nilpotent fundamental group).
\endproclaim

We refer the reader to \cite {AX1} for more precise formulation and proof
of this cobordism theorem.

This result allows one to study complex surfaces and 
Cauchy--Riemannian\break
3-manifolds at their infinity by using decomposition of such 4-manifolds
along homology 3-spheres and applying the gauge theory together with homology
(cobordism) invariants. We mention that, due to Milnor \cite{Ml}, all 
Seifert homology 3-spheres can be seen as the boundaries at infinity of (geometrically finite)
complex hyperbolic 2-orbifolds. Along this line, one can also investigate
the following question (which is also related to Problem 5.4):
 
\proclaim{Problem 4.4}
Are there Cauchy--Riemannian structures on homology 3-spheres of plumbing
type or on real hyperbolic homology 3-spheres? 
\endproclaim   
 
One can be especially interested in this question for homology spheres
obtained by splicing of two Seifert homology spheres along their
singular fibers, see \cite {Sa4}. Another interesting fact (due to Livingston-Myers 
construction \cite{My}) is that any homology 3-sphere is homology
cobordant to a hyperbolic one.
On the other hand,
as it was shown by C.T.C.Wall \cite{Wa}, 
the assignment of the appropriate geometry (when available) gives
a detailed insight into the intrinsic structure of a complex surface.
We mention here Yau's uniformization theorem \cite{Y} which implies that
every smooth complex projective 2-surface $M$ with positive canonical bundle
and satisfying the topological condition that
$\chi(M)=3\operatorname{Signature}(M)$, is a complex hyperbolic manifold.
The necessity of homology sphere decomposition in dimension
four is due to M.Freedman and L.Taylor result  (\cite{ FT}):
\vskip5pt
{\it Let $M$ be a simply connected 4-manifold with intersection
form $q_M$ which decomposes as a direct sum $q_M=q_{M_1}\oplus q_{M_2}$, where
$M_1, M_2$ are smooth manifolds. Then the manifold $M$ can be represented as
a connected sum $M=M_1\#_{\Sigma}M_2$ along a homology sphere $\Sigma$.}
\vskip5pt

We refer to  \cite{Sa1-Sa4}  and \cite{Mat} for recent advances in this direction, in particular,
for results on 
Floer homology of homology 3-spheres and a new Saveliev's (presumably,
homology cobordism) invariant based on Floer homology.

\head 5. Homeomorphisms induced by group isomorphisms
\endhead

The main problem we are concerned in this section is about
geometric realizations of isomorphisms of geometrically finite discrete groups 
$G, H\subset \isom X$ of isometries of a negatively curved space $X$:

\proclaim {Problem 5.1} Given an isomorphism $\varphi\col G\ra H$
of geometrically finite discrete groups 
$G, H\subset \isom X$, find
subsets $X_G, X_H\subset \ov X$ invariant for the action of groups 
$G$ and $H$, respectively,
and an equivariant homeomorphism $f_{\varphi} \col X_G\ra X_H$ which
induces the isomorphism $\varphi $. Determine metric properties of
$f_{\varphi}$, in particular, whether it is either quasisymmetric or
quasiconformal with respect to the given negatively curved metric $d$ in $X$
(or the induced sub-Riemannian structure on the Carnot--Carath\'eodory space 
at infinity $Y=\ov X\bs \{\infty\}$).
\endproclaim

Such type problems were studied by several authors. In the case of
lattices $G$ and $H$ in rank 1 symmetric spaces $X$, G.Mostow \cite{Mo1}
proved in his celebrated rigidity theorem 
that such isomorphisms $\varphi\col G\ra H$ can be extended to inner isomorphisms
of $X$, provided that there is no analytic homomorphism of $X$ onto $PSL(2, \br)$. 
For that proof, it was essential to prove that $\varphi$ can be induced by
a quasiconformal homeomorphism of the sphere at infinity $\ov X$  which is a one point 
compactification of a Carnot group $N$.  Quasiconformal mappings between general Carnot 
groups  have been studied by P.Pansu \cite{P}. For the case of complex
hyperbolic spaces $X=\ch n$ and the Heisenberg groups at their infinity,
foundations of the theory of quasiconformal mappings has been made by
A.Koranyi and M.Reimann \cite {KR1, KR2}.  

The essential component of such
rigidity results is the fact that quasiconformal mappings of Carnot--Carath\'eodory 
spaces are almost everywhere differentiable and preserve the contact structures 
of these spaces (horizontal vector fields). This fact is 
due to Pansu \cite{P} in the case of graded nilpotent groups.
In the case of general Carnot--Carath\'eodory spaces determined by
``horizontal" subbundles of their tangent bundles, this result has been
recently proven by G.Margulis and G.Mostow \cite{MM}. The second 
important ingredient for such rigidity is
the $ACL$-property on lines for quasiconformal mappings, with respect to
the Carnot--Carath\'eodory metric at infinity $\ov X\bs \{\infty\}$, see
\cite{Mo3, KR2, Vo1-2, VG}.

If the groups $G,H\subset \isom X$ are neither lattices nor trivial, the only 
 results on geometric realization of their isomorphisms are known for real
hyperbolic spaces $X$ (of constant negative curvature). In dimension
$\dim X\geq 3$ they are due to P.Tukia \cite {Tu} and use a natural
``type-preserving" condition on the isomorphism $\varphi$, that
parabolic elements are carried out only to parabolic ones. 

In the case of variable negative curvature, it is problematic to use
geometric properties of convex hulls which provide a powerful tool for studying 
metric properties of $G$-equivariant mappings in the real hyperbolic space. 
However, as a first step in solving the above Problem 5.1, we can analyze 
the metric
of $X$ and the induced metrics on horospheres in $X$. Such analysis, 
a generalization
of W.Floyd construction of the group completion \cite{Fl}, and the described above 
analysis of geometrical finiteness in $X$ are crucial components
of our approach to constructing canonical homeomorphisms of subsets in the 
sphere at infinity $\ov X$ , which induce type-preserving isomorphisms $\varphi\col G\ra H$ 
of geometrically finite groups $G$ and $H$ (not necessarily lattices).
In particular, as a first result in this direction, we have the following
isomorphism theorem \cite {A12}:

 \proclaim{Theorem 5.2} Let $\varphi : G\ra H$ be a type preserving isomorphism
 of two non-ele\-men\-ta\-ry geometrically finite discrete subgroups 
$G,H\subset PU(n,1)$. Then
 there exists a unique equivariant ho\-meo\-mor\-phism 
$f_{\varphi}\col \La(G)\ra \La(H)$ of
 their limit sets that induces the isomorphism $\phi$.
 \endproclaim

Upon existence of such equivariant homeomorphism $f_{\varphi}$, the above problem
would  be reduced to the question whether $f_{\varphi}$ is quasisymmetric with
respect to the Carnot--Carath\'eodory metric. In the case of an affirmative answer,
it may be possible to find a global $G$-equivariant homeomorphism (of the sphere at infinity
$\p X$ or even the whole space $\ov X$ inducing the isomorphism $\varphi$.
However, in contrast to the real hyperbolic case,
here we have an interesting fenomenon related to possible noncompactness of the boundary 
$\p M(G)=\om(G)/G$. Namely, even for the simplest case of cyclic groups
$G\cong H \subset PU(n,1)$, the homeomorphic CR-manifolds $\p
M(G)=\sch^n/G$ and $\p M(H)=\sch^n/H$ may be not quasiconformally
equivalent, see \cite {Mn}. So it may be possible to (affirmatively) answer the
following problem:

 \proclaim{Problem 5.3} When are there quasisymmetric homeomorphisms in
a Carnot--Carath\'eodory space $\p X\bs \{\infty\}$ compatible with the action of 
discrete geometrically
finite groups $G,H\subset \isom X$ but quasisymmetrically non-extendable to the whole 
space?
 \endproclaim

We note that, besides the metrical (quasisymmetric) part of this
problem, some topological obstructions for extensions of equivariant
homeomorphisms $f_{\varphi}$ of the limit sets, 
$f_{\varphi}\col \La(G)\ra \La(H)$, may exist. It follows
from the following example.

\proclaim{Example 5.4} Let $G\subset PU(1,1)\subset PU(2,1)$ and 
$H\subset PO(2,1)\subset PU(2,1)$ be two geometrically finite
(loxodromic) groups isomorphic to the fundamental group $\pi_1(S_g)$ of 
a compact oriented surface $S_g$ of genus $g>1$. Then the equivariant
homeomorphism $f_{\varphi}\col \La(G)\ra \La(H)$ cannot be homeomorphically 
extended to the whole sphere $\p \ch 2 \approx S^3$.
 \endproclaim

The obstruction in this example is topological and is due to the fact that 
the quotient
manifolds $\ch 2/G$ and $\ch 2/H$ are not homeomorphic. Namely, these
complex surfaces are disk bundles over the surface $S_g$ and have different Toledo
invariants: $\tau(\ch 2/G)=2g-2$ and $\tau(\ch 2/H)=0$, see \cite {To}.

We remark that some of Carnot--Carath\'eodory spaces are very rigid. In fact, due to
P.Pansu \cite {Pa}, any quasiconformal map on the boundary $\p X$ of
a quaternionic or octanionic hyperbolic space $X$ (which are symmetric
spaces of rank 1) is necessarily an extension of an isometry in $X$. This shows
that any non-trivial homeomorphism $f_{\varphi}\col \La(G)\ra\La(H)$ (non-isometry) 
is non-extendable 
to a quasiconformal map of some open subset of $\p X$. 

In order to attack the above problem, in particular to construct such non-trivial
geometric isomorphisms $\varphi$ in both rigid (like quaternionic and octanionic spaces)
and more flexible (like $\ch n$) spaces, 
one probably may to generalize our block-building method \cite {A1, A5} 
(whose usefulness in
conformal category has been also demonstrated in \cite
{A6, A11}). Another aspect of such
constructions should be 

\proclaim{Problem 5.5} 
Develop geometric combination theorems
generalizing well known Maskit combination theorems for Kleinian groups,
see \cite {Mas, A1}. 
\endproclaim

So far, there are the only simplest versions (free products) of
such combination theorems, see \cite {FG}. The general case of free amalgamated 
product is still unknown even in the complex hyperbolic case, where one can 
however use a week Maskit combination, see \cite{Be}.

As the first steps in investigating Problem 5.1, it appears to be promising 
to study $G$-equivariant homeomorphisms in Carnot groups
$Y=\p X\bs \{\infty\}$, with Carnot--Carath\'eodory metric $\rho$,
 from the following two classes of embeddings 
$f\col A\hookrightarrow Y$, $A\subset Y$. 

The first one consists of well known quasisymmetric embeddings $f$, see
\cite {TV}.

The second class consists of embeddings $f\col A\hookrightarrow Y$, $A\subset Y$,
which generalize the so-called quasi-M\"obius embeddings
introduced in conformal category by V.Aseev \cite {As, AsT} and 
J.V\"ais\"al\"a \cite{V}. They
have bounded distortion of the cross-ratio
(or, equivalently, of the complex cross-ratio in the sense of Koranyi-Reimann
\cite{CR4}), 

$$CR(q)=\rho(x_1,x_2)\rho(x_3,x_4)\rho(x_1,x_3)^{-1}\rho(x_2,x_4)^{-1}
\,,\tag{5.6}$$
of quadruples $q=\{x_1,x_2,x_3,x_4\}\in (Y)^4$ in the following sense. 

Let $\omega\col \br_+\ra\br_+$ be a given homeomorphism,
$\omega(0)=0$. Then, for any quadruple $q\in A^4$ and its $f$-image $f(q)=\{f(x_1,\ldots,
f(x_4\}\in Y^4$, 
we require that 

$$CR(f(q))\leq \omega(CR(q))\,,\tag{5.7}$$
and call such embeddings $f$ {\it quasi-CR
embeddings}. By continuity, the cross ration can trivially be extended
to the one-point compactification $Y\cup \{\infty\}$, so we can consider
embeddings which are not necessarily fixing $\infty$.

In particular, for homeomorphism
$\omega(t)=M\cdot\eta_{\alpha}(t)$ with given constants $M>0$ and 
$\alpha\geq 1$ where

$$\eta_{\alpha}(t)=\cases t^{\alpha}, \quad \text{for $t\geq 1$}\\
                          t^{1/\alpha}, \quad \text{for $0\leq
t<1$}\,,\endcases$$
we call such embeddings $f$ satisfying (5.7) as $(M,\alpha)$-CR-embeddings.

In particular, this approach can be used to study the geometric realization Problem 5.1 in
the case of free geometrically finite groups $G,H\subset \isom X$ whose
limit sets are discontinua. Here we have the following 

\proclaim{Conjecture 5.8} An embedding $f\col A\hookrightarrow \ov Y$,
$A\subset \ov Y$, in a Carnot--Carath\'eodory space $Y$ is quasi-CR if and only if $f$ preserves
the notion of bounded $\mu$-{\it density} of subsets $\Sa\subset A$.  
\endproclaim

Here a set $\Sa$ is called $\mu$-{\it dense} for some $\mu>1$ if any two
points $a,b\in\Sa$ can be connected by a $\mu$-{\it chain} $\{x_i\}$ of points
$x_i\in\Sa$ where $\lim_{i\to -\infty} x_i=a$, $\lim_{i\to +\infty} x_i=b$
and $\ln (CR(a,x_i,x_{i+1},b))\leq \ln \mu$. Naturally, the $\mu$-density 
property is a characteristic property of the limit discontinua for
geometrically finite groups. We also suspect that that property
distinguishes parabolic subgroups of the groups $G$ and $H$ which cannot
be quasiconformally conjugate (we expect that they should have non-maximal rank).
On the other hand, we expect that an embedding $f\col A\hookrightarrow \ov Y$
 of a $\mu$-dense set $A$ in a
Carnot--Carath\'eodory space $Y$ is quasi-CR if and only if it is an $(M,\alpha)$-CR
embedding with some finite $M$ and $\alpha$.
One can expect an affirmative solution of the above problems in the case
of geometrically finite groups whose limit discontinua consist only of
conical points (for the Heisenberg group $Y$, compare \cite {Mn}). For
study of these problems, one can use methods developed by Koranyi and Reimann
\cite{KR1-4}, Vodop'yanov \cite{Vo1,VG} as well as a generalization of
the technique known as Sullivan's microscope (see the next section).

\head 6. Deformations of discrete groups and Sullivan's stability
\endhead         

Here we concern with problems related to deformation spaces of a 
geometrically finite discrete group $G\subset PU(n,1)$ acting by isometries 
in the complex hyperbolic space $\ch n$ and to
 the stability theorem of 
D.Sullivan \cite{Su3} (which has been originally proved for planar Kleinian groups, 
see also \cite{A11, MaG}). We restrict our attention to the complex hyperbolic case 
because other rank 1 spaces with negative variable 
curvature such as quaternionic and octanionic hyperbolic spaces are more
rigid. In particular, due to Corlette's rigidity theorem for harmonic
maps \cite {C2}, finite volume manifolds locally modeled on these spaces are 
super-rigid, completely analogous to Margulis super-rigidity in higher
rank, see \cite{Mg1}. Moreover, due to Gromov and Schoen \cite{GS},
all such finite volume manifolds are arithmetic.

Let $G\subset PU(n,1)$ be a geometrically finite discrete group. Then we
have a fundamental problem concerning
 the deformation space of the inclusion (or of the representation of the group 
$G$ obtained by
restricting a natural inclusion $PU(n,1)\ra PU(n+1,1)$ to $G$). 
The corresponding problem for deformation spaces of a hyperbolically
rigid real hyperbolic lattice $G\subset O(n,1)$, $n\geq 3$, has its
roots in the first construction by the author \cite{A2} of non-
trivial curves in the Teichm\"uller space $T(G)$ of conjugacy classes of
representations of $G$ in $O(n+1,1)$. Then the situation was greatly clarified by
Thurston's ``Mickey Mouse" example \cite{Th} showing that such deformations of a
hyperbolic surface are in fact bendings along geodesics. Since that time, several authors 
studied such deformation spaces, see for example \cite {A4, A6, A8, A9-A11, JM, Ko}. 

In the case of deformations of a complex hyperbolic group $G\subset PU(n,1)$, 
we would like to assume that $G$ is not a uniform lattice, that is the
quotient $\ch n/G$ is noncompact. Otherwise, due to a fundamental result on
local rigidity of deformations \cite {GoM}, the set of ``Fuchsian" representations 
of $G$, $\{\operatorname{Ad} h\col G\ra hGh^{-1},\, h\in PU(n+1,1)\}$,  is a connected component 
$\Cal R_0$
of the variety of representations 
$\Hom (G, PU(n+1,1)$. 

On the other hand, we may deform an isomorphic uniform lattice 
$G'\subset PO(2,1)$, as well as another isomorphic convex co-compact groups 
$G'\subset PU(2,1)$, $G'\cong G$, if such groups have non-elementary subgroups
preserving totally geodesic real planes \cite{AG1}. Such quasiconformal deformations
bend the corresponding complex surfaces along any simple closed geodesic in
its totally real geodesic 2-dim subsurface. For another deformations of groups
$G'\subset PO(2,1)$, see also \cite{GP2}. 
We can also deform the above  groups $G'\subset PO(2,1)$ in $PO(3,1)\subset PU(3,1)$, 
so it makes sense to consider such
deformations in higher dimensions, as it has been done by the author
in the real hyperbolic case, compare \cite{A9-A11}.

The main problem we deal with is as follows.

\proclaim{Conjecture 6.1} Let $G\subset PU(n,1)$ be a geometrically finite
group which is either convex cocompact or having parabolic subgroups of rank at
least 3. Then every its representation $\rho\col G\ra PU(n+k,1)$, $k\geq 0$, that 
close enough to a natural
inclusion, is in fact discrete and faithful and, furthermore, is a quasiconformal 
conjugation. That means that there is an equivariant quasiconformal self-homeomorphism
$f$ of the extended Heisenberg group $\ov {\sch ^{n+k}}$ such that 
$\rho G = f_* G= fGf^{-1}$.
\endproclaim

This conjecture generalizes a remarkable structural stability theorem of
D. Sullivan \cite{Su3} for Kleinian subgroups $G$ in $PSL(2,\bc)$, which shows 
that an algebraic structural stability for holomorphic perturbations
implies a hyperbolicity property for the action of $G$ on its limit set.
In conformal category, we refer to \cite{A11} for our proof of its high dimensional
analog. The complex case provides a possibility to generalize a crucial
Sullivan's argument, the so-called $\la$-Lemma proved by Ma\~ne, Sad and
Sullivan \cite{MSS}. Here topological stability is an intermediate open
problem. We note also that the negative curvature
property is crucial for algebraic stability of convex cocompact
subgroups, due to G.Martin \cite{MaG}:

\proclaim {Theorem 6.2} Let $\{G_t\col t\in\br\}$ be a continuous
deformation of a torsion free convex cocompact group $G\subset \isom X$
in a space $X$ of pinched negative curvature, with convex cocompact
$G_t$. Then $\{G_t\col t\in\br\}$, and its closure in the topology of 
algebraic convergence in $\isom X$, consists entirely of isomorphic groups.
\endproclaim

For geometrically finite groups without parabolics, it is possible to use a generalization 
of the so-called Sullivan's microscope (see \cite{Su2, A11}) to construct
probable quasiconformal conjugation of the group actions on their limit sets. It
is possible because one can find expanding Heisenberg coverings of the limit set
of such a group, which  defines a good Cauchy--Riemannian dynamics there.
As another result important for a generalization of Sullivan's arguments, we mention
an uniformization theorem for Cauchy--Riemannian structures which is due to Falbel 
and Gusevskii \cite{FG}. 

Finally, we would like to mention another problems which are linked with Conjecture 6.1
 and, in dimension $n=2$, with Problem 3.14 about geometrical finiteness.
These problems are related to the boundaries of the deformation spaces 
$\Cal R(G,k)=\Hom (G, PU(n+k,1)$  
and $\Cal T(G,k)=\Cal R(G,k)/PU(n+k,1)$, $k\geq 0$,
to the number of their connected components, discreteness and geometrical finiteness
of boundary representations. In particular:

\proclaim {Problem 6.3} Is any representation $\rho$ obtained as the limit of 
geometrically finite
representations $\rho_i\in \Cal R(G,k)$ discrete and faithful, that is
$\rho(G)\subset PU(n+k)$ is a discrete group isomorphic to a given geometrically 
finite group $G\subset PU(n,1)$?
\endproclaim

We note that, due to G.Martin \cite{MaG}, the subset of discrete faithful
representations is closed in $\Cal R(G,k)$. 

Another question is about the boundary of the Teichm\"uller space 
$\Cal T(G,k)$ of a convex cocompact group $G\subset PU(n,1)$, $n\geq 2,
k\geq 0$. 

\proclaim {Problem 6.4} Are there convex cocompact faithful representations 
$\rho_t\in\Cal R(G,k)$, \break $t\in \br$,
of a convex cocompact group $G\subset PU(n,1)$ which converge to a 
boundary dis\-cre\-te
representation $\rho$ whose image $\rho(G)\subset PU(n+k)$ has accidental parabolic 
ele\-ments?
\endproclaim

We remark that I.Belegradek \cite{Be} constructed a discrete faithful
representation $\rho$ of the fundamental group $G$ 
of a compact hyperbolic surface into $PU(2,1)$ such that $\rho(G)$ has
parabolic elements. This construction uses a Maskit combination theorem,
and it is unclear whether $\rho(G)$ lies in the boundary $\p \Cal R_0(G)$
of a component of the variety of convex co-compact representations $G\to PU(2,1)$.
However the Teichm\"uller space $\Cal T(G)=\Cal T(G,0)$ does have
``cusps" due to a recent result by Apanasov and Gusevskii \cite{AG2}: 

\proclaim {Theorem 6.5} Let $G\subset PO(2,1)\subset PU(2,1)$ be a
uniform lattice isomorphic to the fundamental group of a closed surface
$S_g$ of genus $p\geq 2$. Then, for any simple closed geodesic
$\a\subset S_p=H^2_{\br}/G$, there is a continuous deformation $\rho_t=f^*_t$ 
induced by $G$-equivariant quasiconformal homeomorphisms 
$f_t: \ch 2 \to \ch 2$ whose limit representation $\rho_{\infty}$ corresponds
to a boundary cusp point of the Teichm\"uller space $\Cal T(G)$. In other
words,
the boundary group $\rho_{\infty}(G)$ has an accidental parabolic element 
$\rho_{\infty}(g_{\a})$ where $g_{\a}\in G$ represents the simple closed geodesic 
$\a\subset S_p$.
\endproclaim

We note that, due to our construction of such continuous quasiconformal deformations
in \cite {AG2},
such deformations are independent if the corresponding simple closed geodesics $\a_i\subset S_p$
are disjoint. It implies the existence of a boundary group in $\p \Cal T(G)$ with
``maximal" number of non-conjugate accidental parabolic subgroups:

\proclaim {Corollary 6.6} Let $G\subset PO(2,1)\subset PU(2,1)$ be as in
Theorem 6.5. Then there is a continuous deformation 
$R\,\col\, \br^{2p-2}\to \Cal T(G)$ whose boundary group
$G_{\infty}=R(\infty)(G)$ has $2p-2$ non-conjugate accidental parabolic
subgroups.
\endproclaim

Finally, we mention another aspect of the intrigue Problem 3.14:

\proclaim {Problem 6.7} Construct a geometrically infinite (finitely
generated) group $G\subset PU(n,1)$, $n\geq 2$, whose limit set is the whole sphere
at infinity, $\La(G)=\p \ch n=\ov {\sch^n}$, and which is the limit of
convex cocompact groups $G_i\subset PU(n,1)$ from the Teichm\"uller space 
$\Cal T(G)$ of a convex cocompact group $G\subset PU(n,1)$. Is that
possible for a Schottky group $G$?
\endproclaim

\vfil
\eject

\def\ref#1{[#1]}
\eightpoint
\parindent=36pt

\head  REFERENCES\endhead
\bigskip

\frenchspacing

\item{\ref{Ah}} Lars V.Ahlfors, Fundamental polyhedra and limit point sets
of Kleinian groups. - Proc. Nat. Acad. Sci. USA, {\bf 55}(1966), 251-254.

\item{\ref{As}} V.V.Aseev, Quasi-symmetric embeddings and mappings
with bounded distortion of mo\-du\-li. - Sib. Matem. Zh./VINITI, Dep 7190-84, 1984
(Russian).

\item{\ref{AsT}} V.V. Aseev and D.A. Trotsenko, Quasi-symmetric embeddings,
quadruples of points, and mo\-du\-li distortion. - Sibirsk. Matem. Zh. {\bf
28:4} (1987), 32-38 (Russian); Engl. Transl.: Siberian Math. J. {\bf 28}
(1987).

\item{\ref{Au}} L.~Auslander, Bieberbach's theorem on space groups and discrete 
uniform subgroups of Lie groups, II. - Amer. J. Math., {\bf 83} (1961), 276-280.

\item{\ref{A1}}  Boris Apanasov, Discrete groups in Space and
Uniformization  Problems. - Math. and Appl. {\bf 40}, Kluwer
Academic Publishers, Dordrecht, 1991.

\item{\ref{A2}}  \underbar{\phantom{Apanasov}}, Nontriviality  of Teichm\"uller space for
Kleinian  group in  space. -  Riemann Surfaces and Related
Topics: Proc. 1978 Stony  Brook Conference (I.~Kra and
B.~Maskit, Eds), Ann. of  Math. Studies {\bf 97},
Princeton Univ. Press, 1981, 21-31.

\item{\ref{A3}}  \underbar{\phantom{Apanasov}}, Geometrically finite hyperbolic
structures on manifolds. - Ann. of Glob. Ana\-ly\-sis and Geom. {\bf 1:3} (1983),
1-22.

\item{\ref{A4}}    \underbar{\phantom{Apanasov}},   Thurston's bends and geometric
deformations of conformal structures. - Comp\-lex Analysis and
Applications'85,  Publ. Bulgarian Acad. Sci.,
Sofia, 1986, 14-28.

\item{\ref{A5}}  \underbar{\phantom{Apanasov}}, Quasisymmetric embeddings of a closed ball
inextensible in neighborhoods of any boundary point. - Ann. Acad. Sci. Fenn. Ser. A I 
Math. {\bf 14} (1989), 243-255.

\item{\ref{A6}}  \underbar{\phantom{Apanasov}},  Nonstandard  uniformized  conformal
structures on hyperbolic manifolds. - \break Invent. Math. {\bf 105} (1991), 137-152.

\item{\ref{A7}}  \underbar{\phantom{Apanasov}}, Kobayashi
conformal metric on manifolds, Chern-Simons and $\eta$-invariants. - 
Intern. J. Math. {\bf 2} (1991), 361-382.

\item{\ref{A8}}  \underbar{\phantom{Apanasov}}, Deformations of conformal structures on
hyperbolic manifolds. - J. Diff. Geom. {\bf 35}
(1992), 1-20.

\item{\ref{A9}}   \underbar{\phantom{Apanasov}}, Quasisymmetric knots
and Teichm\"uller spaces.- Russian Acad. Sci. Dokl. Math. {\bf 343}
(1995).

\item{\ref{A10}}   \underbar{\phantom{Apanasov}}, Varieties of discrete representations of
hyperbolic 3-lattices.- 
M.I.Kargopolov Memorial Volume (Yu.L.Ershov a.o., eds), W. de Gruyter-Verlag, Berlin,
1996, 7-20.

\item{\ref{A11}} \underbar{\phantom{Apanasov}}, Conformal geometry of
discrete groups and manifolds. - W. de Gruyter, Berlin - New York, 
to appear.

\item{\ref{A12}} \underbar{\phantom{Apanasov}}, Canonical homeomorphisms in
Heisenberg groups induced by isomorphisms of discrete subgroups of $PU(n,1)$.
 - Russian Acad. Sci. Dokl. Math., to appear.

\item{\ref{AG1}}  Boris Apanasov and Nikolay Gusevskii,  Bending
deformations of complex hyperbolic surfaces. - In pre\-pa\-ra\-tion.

\item{\ref{AG2}}  \underbar{\phantom{Apanasov}},  The boundary of
Teichm\"uller space of complex hyperbolic surfaces. - In pre\-pa\-ra\-tion.

\item{\ref{AX1}}  B.Apanasov and X.Xie, Geometrically finite 
complex hyperbolic manifolds. - Preprint, Univ. of Oklahoma, 1995, 51 pp.

\item{\ref{AX2}} \underbar{\phantom{Apanasov}}, Manifolds of negative
curvature and nilpotent groups. - Prep\-rint, Univ. of Oklahoma, 1995.

\item{\ref{AX3}} \underbar{\phantom{Apanasov}}, Discrete isometry groups
of nilpotent Lie groups. - Preprint, Univ. of Oklahoma, 1995.

\item{\ref{AF}} C.S.~Aravinda and F.T.~Farrell,  Rank 1  aspherical manifolds
which do not support any nonpositively curved metric.  -
Comm. Anal. Geom. {\bf 2} (1994), 65-78.

\item{\ref{BGS}} W.~Ballman, M.~Gromov and V.~Schroeder,  Manifolds of
nonpositive curvature. - Birk\-h\"au\-ser, 1985.

\item{\ref{BM}} A. Beardon and B. Maskit, Limit points of Kleinian
groups and finite sided fundamen\-tal polyhedra. - Acta Math. {\bf 132} (1974),
1-12.

\item{\ref{Be}} I.~Belegradek, Discrete surface groups actions with
accidental parabolics on comp\-lex hyperbolic plane. - Preprint, 1995.

\item{\ref{Bo}} B.~Bowditch,
Geometrical finiteness with variable negative curvature. -
Duke J. Math. {\bf 77} (1995), 229-274.

\item{\ref{BuM}} D.~Burns and R.~Mazzeo, On geometry of cusps for
$SU(n,1)$. - Preprint, Univ. of Michigan, 1994.

\item{\ref{BS}} D.~Burns and S.~Shnider, Spherical hypersurfaces in complex
manifolds. - Inv. Math., {\bf 33}(1976), 223-246.

\item{\ref{C1}} K. Corlette, Hausdorff dimensions of limit sets. -
Invent. Math. {\bf 102} (1990), 521-542.

\item{\ref{C2}}  \underbar{\phantom{Apanasov}}, Archimedian
superrigidity and hyperbolic geometry. - Ann. of Math. {\bf 135} (1992),
165-182.

\item{\ref{CG}} S. Chen and L. Greenberg, Hyperbolic spaces. -
Contributions to Analysis, Aca\-de\-mic Press, New York, 1974, 49-87.

\item{\ref{DS}} S.K. Donaldson and D. Sullivan, Quasiconformal 
4-manifolds. -  Acta Math. {\bf 163} (1989), 181-252.

\item{\ref{EHS}} P.Eberlein, U.Hamenst\"adt and V.Schroeder, Manifolds
of nonpositive curvature. - Proc. Symp. Pure Math. {\bf 54:3}, 1993,
179-227.

\item{\ref{EMM}} C. Epstein, R. Melrose and G. Mendoza, Resolvent of the Laplacian
on strict\-ly pseu\-do\-con\-vex domains. - Acta Math. {\bf 167} (1991), 1-106.

\item{\ref{FG}} E. Falbel and N. Gusevskii, Spherical CR-manifolds of
dimension 3. - Bol. Soc. Bras. Mat. {\bf 25} (1994), 31-56.

\item{\ref{FT}}  M.H. Freedman and  L. Taylor, $\Lambda$-splitting 
4-manifolds. -
Topology {\bf 16} (1977), 181-184.

\item{[GM]} F.W. Gehring and G.J. Martin, Discrete quasiconformal groups. -
Proc. London Math. Soc. {\bf 55} (1987), 331-358.

\item{\ref{Go}} W. Goldman, Complex hyperbolic geometry. - 
Oxford Univ. Press, to appear. 

\item{\ref{GoM}} W. Goldman and J. Millson, Local rigidity of discrete
groups acting on complex hyper\-bo\-lic space. - Invent. Math., {\bf 88}
(1987), 495-520.

\item{\ref{GP1}}  W. Goldman and J. Parker, Dirichlet polyhedron for
dihedral groups acting on comp\-lex hyperbolic space. -
J. of Geom. Analysis {\bf 2:6} (1992), 517-554.

\item{\ref{GP2}} \underbar{\phantom{Apanasov}}, Complex hyperbolic ideal
triangle groups.- J. reine angew. Math. {\bf 425} (1992), 71-86.

\item{\ref{Gr}} M. Gromov,  Carnot--Carath\'eodory spaces seen
 from within.- Preprint IHES, Bures-sur-Yvette, 1994.

\item{\ref{JM}}  D. Johnson and  J. Millson,  Deformation  spaces
associated  to  compact  hyperbolic  ma\-ni\-folds. - Discrete
Groups in Geometry  and Analysis,%:  Papers in  Honor of G.D.
%Mostow on  His Sixtieth Birthday,  Ed. R. Howe, Birkhauser,
Boston, 1987, 48-106.

\item{\ref{KR1}} A. Koranyi and M. Reimann, Quasiconformal mappings on the
Heisenberg group. - Invent. Math. {\bf 80} (1985), 309-338.

\item{\ref{KR2}} \underbar{\phantom{Apanasov}}, Foundations for the
theory of quasiconformal mappings on the Hei\-sen\-berg group. - Adv. Math.
{\bf 111} (1995), 1-87.

\item{\ref{KR3}} \underbar{\phantom{Apanasov}}, Contact transformations
as limits of symplectomorphisms. - C. R. Acad. Sci. Paris {\bf 318}
(1994), 1119-1124.

\item{\ref{KR4}} \underbar{\phantom{Apanasov}}, The complex cross ratio
on the Heisenberg group.- L'Ens. Math. {\bf 33} (1987), 291-300.

\item{\ref{Ko}} C. Kourouniotis,
Deformations of hyperbolic structures.- Math. Proc. Cambr. Phil. Soc.
 {\bf 98} (1985), 247-261.

\item{\ref{LB}} Claude LeBrun, Einstein metrics and Mostow rigidity. -
Preprint, Stony Brook, 1994.

\item{\ref{MSS}} R.Ma\~ne, P.Sad and D.Sullivan,
On the dynamics of rational maps.- Ann. Sci. \'Eco\-le Norm. Sup. {\bf
16} (1983), 193-217.

\item{\ref{Md}} A. Marden, The geometry of finitely generated Kleinian
groups. - Ann. of Math. {\bf 99} (1974), 383-462.

\item{\ref{Mg1}} G.A. Margulis,
Discrete groups of motions of manifolds of nonpositive curvature. -
Amer. Math. Soc. Transl. {\bf 109} (1977), 33-45.

\item{\ref{Mg2}} \underbar{\phantom{Apanasov}}, Free properly discontinuous
groups of affine transformations. - Dokl. Acad. Sci. USSR {\bf 272} (1983),
937-940.

\item{\ref{MM}} G.A. Margulis and G.D. Mostow, The differential of a
quasiconformal mapping of a Car\-not-Ca\-ra\-th\'eodory space. - Geom. Funct.
Anal. {\bf 5} (1995), 402-433.

\item{\ref{MaG}} G. Martin, On discrete isometry groups of negative
curvature.- Pacif. J. Math. {\bf 160} (1993), 109-127.

\item{\ref{Mas}} B.Maskit, Kleinian groups. - Springer-Verlag, 1987.

\item{\ref{Mat}} R. Matveyev, A decomposition of smooth simply-connected
$h$-cobordant 4-manifolds. - Preprint, Michigan State Univ. at
E.Lansing, 1995.

\item{\ref{Ml}} J. Milnor, On the 3-dimensional Brieskorn manifolds
$M(p,q,r)$.- Knots, groups and \break 3-manifolds, Ann. of Math.
Studies {\bf 84}, Princeton Univ. Press, 1975, 175--225.

\item{\ref{Mn}} R. Miner, Quasiconformal equivalence of spherical CR manifolds.-
Ann. Acad. Sci. Fenn. Ser. A I Math. {\bf 19} (1994), 83-93.

\item{\ref{Mo1}} George D. Mostow, Strong rigidity of locally symmetric
spaces.- Princeton Univ. Press, 1973.

\item{\ref{Mo2}} \underbar{\phantom{Apanasov}},  On a remarkable class of polyhedra in
complex hyperbolic space. - Pacific J. Math., {\bf 86} (1980), 171-276.

\item{\ref{Mo3}} \underbar{\phantom{Apanasov}}, A remark on
quasiconformal mappings on Carnot groups.- Mich. Math. J. {\bf 41}
(1994), 31-37.

\item{\ref{My}} R. Myers, Homology cobordisms, link concordances, and
hyperbolic 3-manifolds.- Trans. Amer. Math. Soc.  {\bf 278} (1983), 271-288.

\item{\ref{Pa}} P. Pansu, M\'etriques de Carnot--Carath\'eodory et
quasiisom\'etries des espaces sym\-m\'et\-ries de rang un.- Ann. Math. {\bf
129} (1989), 1-60.

\item{\ref{Ph}} M.B. Phillips, Dirichlet polyhedra for cyclic groups in complex
hyperbolic space. - Proc. Amer. Math. Soc. {\bf 115} (1992), 221-228.

\item{\ref{Sa1}}  Nikolai Saveliev, Floer homology and 3-manifold invariants.-
Ph.D. Thesis, Univ. of Ok\-la\-ho\-ma, 1995.

\item{\ref{Sa2}}  \underbar{\phantom{Apanasov}}, Floer homology and invariants of 
homology cobordism.- Preprint, Univ. of Oklahoma, 1995.

\item{\ref{Sa3}}  \underbar{\phantom{Apanasov}}, On homology cobordisms of plumbed 
homology 3-spheres.- Preprint, Univ. of Oklahoma, 1995.

\item{\ref{Sa4}}  \underbar{\phantom{Apanasov}}, Notes on homology cobordisms of 
plumbed homology 3-spheres.- Preprint, Univ. of Michigan, 1995.  

\item{\ref{Su1}} D. Sullivan, Hyperbolic geometry and
homeomorphisms.- Geometric topology, Acad. Press, 1979.

\item{\ref{Su2}}  \underbar{\phantom{Apanasov}}, Seminar  on  conformal  and 
hyperbolic geometry.- Preprint  IHES, 1982.

\item{\ref{Su3}} \underbar{\phantom{Apanasov}}, Quasiconformal homeomorphisms and 
dynamics, II: Structural stability implies hyperbolicity for Kleinian 
groups. - 
Acta Math. {\bf 155} (1985), 243--260.

\item{\ref{Te}} A.V. Tetenov, The discontinuity set for
a Kleinian group and topology of its Klei\-nian manifold. -
Intern. J. Math. {\bf 4} (1993), 167-177.

\item{\ref{Th}} W.P. Thurston, The geometry and topology of
three-manifolds. - Lect. Notes, Prince\-ton Univ., 1979.

\item{\ref{To}} D. Toledo, Representations of surface groups on complex
hyperbolic space. - J. Diff. Geom. {\bf 29} (1989), 125-133.

\item{\ref{Tu}} P. Tukia, On isomorphisms of geometrically finite
Kleinian groups. - Publ. Math. IHES {\bf 61} (1985), 171-214.

\item{\ref{TV}} P. Tukia and J. V\"ais\"al\"a, Quasisymmetric embeddings
of metric spaces.- Ann. Acad. Sci. Fenn. Ser. A I Math. {\bf 5} (1980),
97-114.

\item{\ref{V}}  J. V\"ais\"al\"a, Quasim\"obius maps.- J. Analyse Math.
{\bf 50} (1984/85), 218-234.

\item{\ref{Vo1}} S.K. Vodop'yanov, Quasiconformal mappings on Carnot
groups.- Russian Dokl. Math., to appear.

\item{\ref{Vo2}} \underbar{\phantom{Apanasov}}, Monotone functions on 
Carnot groups.-  Siberian Math. J. {\bf 37:5} (1996), to appear.

\item{\ref{VG}} S.K. Vodop'yanov and A.V.Greshnov, Analytic properties of 
quasiconformal mappings on Carnot groups. - Siberian Math. J. {\bf 36:6}
(1995), 1317-1327.

\item{\ref{VU}} S.K. Vodop'yanov and A.D. Ukhlov, Weakly contact transformations  
and a chan\-ge-of-variable formula on nilpotent groups. - 
Russian Acad. Sci. Dokl. Math. {\bf 341:4} (1995), 439-441.

\item{\ref{Wa}} C.T.C.~Wall, Geometric structures on compact complex analytic surfaces.
- Topology, {\bf 25} (1986), 119-153. 

\item{\ref{Wo}} J. Wolf, Spaces of constant curvature.- Publ. or Perish,
Berkeley, 1977.

\item{\ref{Y}} S.T.Yau, Calabi's conjecture and some new results in
algebraic geometry.- Proc. Nat. Acad. Sci {\bf 74} (1977), 1798-1799.

\enddocument